\documentclass[12pt]{article}
\usepackage{amsmath}
\usepackage{amssymb}

\newcommand{\eop}{\bigstar}  % end-of-proof

\newenvironment{Proof}{\noindent{\bf Proof.}}{\par\bigskip} 

\newtheorem{THEOREM}{Theorem}[section]

\newtheorem{Conclusion}[THEOREM]{Conclusion}

\newtheorem{Hypothesis}[THEOREM]{Hypothesis}

\newtheorem{LEMMA}[THEOREM]{Lemma}

\newtheorem{Main Theorem}[THEOREM]{Main Theorem}
\newenvironment{main Theorem}{\begin{Main Theorem}} 
{\end{Main Theorem}}

\newtheorem{Theorem}[THEOREM]{Theorem}

\newtheorem{Definition}[THEOREM]{Definition}

\newtheorem{Conventions}[THEOREM]{Conventions}

\newtheorem{Main Definition}[THEOREM]{Main Definition}
\newenvironment{main definition}{\begin{Main Definition}}
{\end{Main Definition}}

\newtheorem{Lemma}[THEOREM]{Lemma}

\newtheorem{Notation}[THEOREM]{Notation}

\newtheorem{Convention}[THEOREM]{Convention}

\newtheorem{Note}[THEOREM]{Note}

\newtheorem{Observation}[THEOREM]{Observation}

\newtheorem{Remark}[THEOREM]{Remark}

\newtheorem{Main Fact}[THEOREM]{Main Fact}
\newenvironment{main Fact}{\begin{Main Fact}}{\end{Main Fact}}

\newtheorem{Fact}[THEOREM]{Fact}

\newtheorem{Subfact}[THEOREM]{Subfact}

\newtheorem{Claim}[THEOREM]{Claim}

\newtheorem{Main Claim}[THEOREM]{Main Claim}
\newenvironment{main claim}{\begin{Main Claim}}{\end{Main Claim}}

\newtheorem{Crucial Claim}[THEOREM]{Crucial Claim}
\newenvironment{crucial claim}{\begin{Crucial Claim}}{\end{Crucial Claim}}

\newtheorem{Subclaim}[THEOREM]{Subclaim}

\newtheorem{Corollary}[THEOREM]{Corollary}

\newtheorem{Example}[THEOREM]{Example}

\newtheorem{Problem}[THEOREM]{Problem}

\newtheorem{Proposition}[THEOREM]{Proposition}

\newtheorem{Discussion}[THEOREM]{Discussion}

\newenvironment{Proof of the Subfact}
{\noindent{\bf Proof of the Subfact.}}{\par\bigskip}

\newenvironment{Proof of the Theorem}
{\noindent{\bf Proof of the Theorem.}}{\par\bigskip}

\newenvironment{Proof of the Proposition}
{\noindent{\bf Proof of the Proposition.}}{\par\bigskip}

\newenvironment{Proof of the Conclusion}
{\noindent{\bf Proof of the Conclusion.}}{\par\bigskip}

\newenvironment{Proof of the Observation}
{\noindent{\bf Proof of the Observation.}}{\par\bigskip}

\newenvironment{Proof of the Fact}
{\noindent{\bf Proof of the Fact.}}{\par\bigskip}

\newenvironment{Proof of the Lemma}
{\noindent{\bf Proof of the Lemma.}}{\par\bigskip}

\newenvironment{Proof of the Claim}
{\noindent{\bf Proof of the Claim.}}{\par\bigskip}

\newenvironment{Proof of the Corollary}
{\noindent{\bf Proof of the Corollary.}}{\par\bigskip}

\newenvironment{Proof of the Subclaim}
{\noindent{\bf Proof of the Subclaim.}}{\par\medskip}

\newenvironment{Proof of the Main Claim}
{\noindent{\bf Proof of the Main Claim.}}{\par\bigskip}

\newenvironment{Proof of the Crucial Claim}
{\noindent{\bf Proof of the Crucial Claim.}}{\par\bigskip}

\catcode`\@=11
\def\@begintheorem#1#2{\rm \trivlist \item[\hskip \labelsep{\bf #1\ #2.}]}
\def\@opargbegintheorem#1#2#3{\rm \trivlist
      \item[\hskip \labelsep{\bf #1\ #2\ (#3).}]}
\catcode`\@=12

\catcode`\@=11
\def\@&{\hskip2pt \&\hskip2pt}
\catcode`\@=12

\newcommand{\into}{\rightarrow}

\newcommand{\PP}{{\cal P}}

\newcommand{\concat}{\kern-.25pt\raise4pt\hbox{$\frown$}\kern-.25pt}

\newcount\skewfactor
\def\mathunderaccent#1#2 {\let\theaccent#1\skewfactor#2
\mathpalette\putaccentunder}
\def\putaccentunder#1#2{\oalign{$#1#2$\crcr\hidewidth
\vbox to.2ex{\hbox{$#1\skew\skewfactor\theaccent{}$}\vss}\hidewidth}}
\begin{document}

\title{An application of CAT}

\author{Mirna D\v zamonja and Jean Larson}

\maketitle

\begin{abstract}
We comment on a question of Justin Moore on colourings of pairs
of nodes in an Aronszajn tree.
\end{abstract}

%\section{}

\begin{Lemma}\label{lem:chain}
For any uncountable antichain $A$ in an Aronszajn tree $T$,
there is an infinite chain $C$ in $T$ such that every element
of $C$ is the meet of two elements of $A$.
\end{Lemma}

\begin{Proof} Define a sequence $\langle (t_i, u_i, v_i, B_i):\,i<\omega\rangle$
by recursion as follows. To start let $t_0$ be an arbitrary
element of $A$ (note that $t_0$ cannot be the root of $T$). Consider all 
$t_0\wedge a$ for $a\in A$ and notice that this is a countable set.
Hence there is $u_0<_T t_0$ such that $\{a\in A:\,u_0= a\wedge t_0\}$ is uncountable.
Note that $u_0$ has at least two distinct immediate successors in $T$, so let
$v_0$ be an immediate successor of $u_0$ which is incompatible with $t_0$ and which
has uncountably many extensions in $A$. Denote the set of extensions of $v_0$ in $A$
by $B_0$.

At the stage $i=j+1$ we choose $t_i\in B_j$ and $u_i<_T t_i$ such that 
$\{a\in B_j:\,u_i= a\wedge t_i\}$ is uncountable. Choose $v_0$ to be an immediate successor
of $u_i$ which is incompatible with $t_i$ and which
has uncountably many extensions in $B_j$. Let this set of extensions be $B_i$. Note that
$u_j<_T v_j <_T u_i <_T t_i$.

At the end the set $\{u_i:\,i<\omega\}$ is an infinite chain such that $u_i= t_i\wedge t_{i+1}$.
$\eop_{\ref{lem:chain}}$
\end{Proof}

\begin{Definition} (1) A {\em subtree} of an $\omega_1$-tree $T$ will mean an uncountable
meet closed subset of $T$.

{\noindent (2)} A subtree $S$ of an $\omega_1$-tree is {\em binary} if every node of $S$
has at most two distinct immediate successors in $S$.

{\noindent (3)} An $\omega_1$-tree $T$ is {\em binarisable} if every subtree of $T$ has
a binary subtree.

{\noindent (4)} For a tree $T$ we let $T^{[2]}=\{\{s,t\}:\,s<_T t\}$.
\end{Definition}

The above definition will be used in the context of Aronszajn trees, so the case of trivial
binary trees, namely uncountable branches, will be avoided. We will use the following
statement introduced in \cite{AbrSh} and used in \cite{Stevo}:

{\bf Colouring Axiom for Trees (CAT):} For any partition $T=K_0\cup K_1$ of
an Aronszajn tree $T$, there is an uncountable set $X\subseteq T$ and $i<2$
such that $x\wedge y\in K_i$ for all distinct $x, y\in X$.

We remark that by repeated applications of CAT one obtains that for any partition of
an Aronszajn tree $T$ into finitely many pieces there is an uncountable set $X\subseteq T$
such that $x\wedge y$ lie in the same piece of the partition, for all distinct $x, y\in X$.

\begin{Theorem}\label{thm:p} (CAT) For every binarisable special Aronszajn tree
$T\subseteq {}^{\omega_1>} p$ there is a colouring $c$ of $T^{[2]}$ into $p+1$
colours such that every subtree $S$ of $T$ realises at least 3 colours.
\end{Theorem}

\begin{Proof} Let $T\subseteq {}^{\omega_1>} p$ be a given binarisable special Aronszajn tree.
Let $T=\bigcup_{n<\omega} A_n$ witness that $T$ is special, so each $A_n$
is an antichain. We shall assume that $A_n$'s are disjoint, and for $t\in T$ let
$n(t)$ be $n$ such that $t\in A_n$. Note that if $s<_T t$ then $n(s)\neq n(t)$.
Define $c:\, T^{[2]} \to p+1$ for $s<_T t$ by:
\[
c(\{s,t\})=\begin{cases}
p & \text{if } n(s)<n(t),\\
i & \text{if } n(s)>n(t)\text{ and }t(\lg(s))=i.\\
\end{cases}
\]
Here $\lg(s)$ is the order type of the domain of $s$. Let now $S$ be any subtree of $T$. Let
$S'\subseteq S$ be a binary subtree of $S$. Order the subsets of $p$ by $<^\ast$. Let
$d:\,S'\into \PP(p)$ be given by 
\[
d(s)=\{j:\,s\concat j\mbox{ has an extension in }S'\}.
\]
Applying CAT we obtain an uncountable $X\subseteq S'$ and $J\subseteq p$ such that 
for all $x\neq y\in X$ we have $d(x\wedge y)=J$. Since $S'$ is binary
the cardinality of $J$ is at most 2, but since $X$ is uncountable and
$T$ has no uncountable branches, $J$ must have cardinality exactly 2.
Let $S''$ be the set of meets
of distinct elements of $X$ and note that this is a subtree of $S'$, hence of $S$.
Also, every node of $S''$ has the property that $s\concat j$
has an extension in $X$ iff $j\in J$. 

To complete the proof we shall show that for all $j\in J\cup\{p\}$ there is
$\{s,t\}\in S''^{[2]}$ such that $c(\{s,t\})=j$.
Let $m$ be such that $S''\cap A_m$ is uncountable and apply Lemma \ref{lem:chain} to
$S''\cap A_m$ to obtain an infinite chain $C\subseteq S''$ such that 
\[
(\forall x\in C)(\exists y,z\in S''\cap A_m)\,x=y\wedge z.
\]
Note that $s\neq t\in C$ implies $n(s)\neq n(t)$, hence there are $s_p$ and $t_p$ in $C$
such that $c(\{s_p,t_p\})=p$. Let $x\in C$ be such that $n(x)>m$. Find $y\neq z\in
S''\cap A_m$ such that $x=y\wedge z$. Since $x\in S''$, $d(x)=J$. Hence
$\{c(\{x,y\}), c(\{x,z\})\}=J$.
$\eop_{\ref{thm:p}}$
\end{Proof}

With $p=2$ Theorem \ref{thm:p} states that under CAT,
for every special Aronszajn subtree of $\subseteq {}^{\omega_1>} 2$, there is a colouring
of $T^{[2]}$ into $3$ colours such that every subtree $S$ of $T$ realises all 3 colours.
Since CAT is implied by PFA, which also implies that all Aronszajn trees are special,
this answers negatively Question 9.6. of \cite{justin}. For further references and results on CAT
see \cite{justin},
especially Theorem 5.2.

\end{document}